\documentclass[a4paper,12pt]{article}
\usepackage{amsfonts,amssymb,amsmath,amsthm}
\usepackage[dvips]{graphics}
\usepackage[dvips]{graphicx}
\usepackage[dvips]{color}

\newtheorem{st}{Statement}[section]

\newtheorem{propo}[st]{Proposition}

\newtheorem{cor}[st]{Corollary}
\newtheorem{cl}[st]{Cleim}
\newtheorem{rem}[st]{Remark}
\newtheorem{thm}[st]{Theorem}
\newtheorem{lemm}[st]{Lemma}

\def\max{{\rm max\,}}
\def\min{{\rm min\,}}
\def\sup{{\rm sup\,}}

\def\card{{\rm card\,}}

\def\ss{{{_{\rm Shift}}\,}}

\def\cP{{\cal P}}

\def\card{{\rm card\,}}
\def\lim{{\rm lim\,}}

\def\dom{{\rm dom\,}}

\def\rng{{\rm rng\,}}

\def\SC{{\rm SC\,}}

\def\Proof:{ \vspace{-1.5mm} {\noindent\it Proof.}}
\def\gm{\vspace{0mm}}
\def\gd{\vspace{0mm}}

\def\Box{\rule{1.5mm}{1.5mm}}

\begin{document}

\renewcommand{\thefootnote}{\fnsymbol{footnote}}

 \title{\gd \gd Rudin-Kisler ordering on the P-hierarchy}
 \author{Andrzej Starosolski}
 \date{\today}
\maketitle

\gd
\begin{abstract}
In \cite{MERudin} M. E. Rudin proved (under CH) that for each P-point $u$ there is a P-point $v$ such that $v>_{RK}u$. In \cite{Blass1} A. Blass improved that theorem assuming MA$^{1)}$ in the place of CH, in that paper he also proved that under 
MA$^{1)}$ \footnotetext{ 1) The theorem was stated under MA, but in fact the Blass proof works also under $\mathfrak{p}=\mathfrak{c}$, which was mentioned in \cite{Blass1}}each RK-increasing sequence of P-points is upper bounded by a P-point. We improve Blass results simultaneously  in 3 directions - we prove it for each class of index $\geq 2$ of P-hierarchy (P-points coincidence with a class $P_2$ of P-hierarchy), assuming $\mathfrak{b}=\mathfrak{c}$ 
in the place of MA and we show that there are at least $\mathfrak{b}$ many Rudin-Kisler incomparable such upper bounds.

\end{abstract}


\gd
\section{Introduction}

We proved in \cite{Star-P-hier} that a class of P-points is precisely a class $\cP_2$ of P-hierarchy 
which is a classification of ultrafilters on $\omega$ into $\omega_1$ disjoint classes. It is natural to ask which properties of the class of P-points are (or are not) also properties of other classes of P-hierarchy.
We have started this work in earlier papers \cite{Star-P-hier} and \cite{Star-Ord-V-P} where also Rudin-Kisler ordering was examined. Here, inspired by papers of M. E. Rudin and A. Blass we continue our investigation. The P-hierarchy is defined by monotone sequential contours, and since this ideas are not widely known here we recall all necessary informations.

In \cite{DolMyn} S. Dolecki and F. Mynard introduced monotone
sequential cascades - special kind of trees - as a tool to describe
topological sequential spaces. Cascades and their contours appeared
to be also an useful tool to investigate certain types of ultrafilters
on $\omega$, namely ordinal ultrafilters and the P-hierarchy
(see \cite{Star-P-hier}, \cite{Star-Ord-V-P}).

The {\it cascade} is a tree $V$, ordered by "$\sqsubseteq $",
 without infinite branches and with the minimal element $\emptyset _V$.
A cascade is $\it sequential$ if for
each non-maximal element of $V$ ($v \in V \setminus \max V$) the set
$v^{+V}$ of immediate successors of $v$ (in $V$) is countably
infinite. We write $v^+$ instead of $v^{+W}$ if it is known in which
cascade the successors of $v$ are considered. If $v \in V \setminus
\max V$, then the set $v^+$ (if infinite) may be endowed with an
order of the type $\omega$, and then by $(v_n)_{n \in \omega}$ we
denote the sequence of elements of $v^+$, and by $v^{(n)W}$ - the
$n$-th element of $v^{+W}$.

The {\it rank} of $v \in V$ ($r_V(v)$ or $r(v)$) is defined
inductively as follows: $r(v)=0$ if $v \in \max V$, and otherwise
$r(v)$ is the least ordinal greater than the ranks of all immediate
successors of $v$. The rank $r(V)$ of the cascade $V$ is, by
definition, the rank of $\emptyset_V$. If it is possible to order
all sets $v^+$ (for $v \in V \setminus \max V$)  so that for each $v
\in V \setminus \max V$ the sequence $(r(v^{(n)})_{n<\omega})$ is
non-decreasing (in other words if for each 
$v\in V\setminus \emptyset_V$ the set $\{v\in (w)^+
: r(v)<\alpha\}$ is finite for each $\alpha<r(w)$),
then the cascade $V$ is {\it monotone}, and we fix
such an order on $V$ without indication. 

For $v\in V$ by $v^\uparrow$ we understand $\{w\in V: v \sqsubseteq w\}$ with preserved order,
if $V$ is a monotone sequential cascade and $U\# \int V$ then by $V^{\downarrow V}$
we understand the biggest monotone sub-cascade of cascade $V$ such that for each element $w\in V^{\downarrow U}$
we have $\max(w^\uparrow) \in U$.

Let $W$ be a cascade, and let $\{V_w: w \in \max W\}$ be a set of
pairwise disjoint cascades such that $V_w \cap W = \emptyset$ for
all $w \in \max W$. Then, the {\it confluence} of cascades $V_w$
with respect to the cascade $W$ (we write $W \looparrowleft V_w$) is
defined as a cascade constructed by an identification of $w \in
\max W$ with $\emptyset_{V_w}$ and according to the following rules:
$\emptyset_W = \emptyset_{W \looparrowleft V_w}$; if $w\in W
\setminus \max W$, then $w^{+ W \looparrowleft V_w} = w^{+W}$; if $w
\in V_{w_0}$ (for a certain $w_0 \in \max W$), then $w^{+ W
\looparrowleft V_w} = w^{+V_{w_0}}$; in each case we also assume
that the order on the set of successors remains unchanged. By $(n)
\looparrowleft V_n$ we denote $W \looparrowleft V_w$ if $W$ is a
sequential cascade of rank 1.

If $\mathbb{U}= \{{u}_s: s \in S \}$ is a family of filters on $X$
and if $p$ is a filter on $S$, then the {\it contour of $\{ {u}_s
\}$ along $p$} is defined by
$$\int_{p} \mathbb{U} = \int_{p}{u}_s =
\bigcup_{P \in p} \bigcap_{s \in P} {u}_s.$$

Such a construction has been used by many authors (\cite{Fro},
\cite{Gri1}, \cite{Gri2}) and is also known as a sum (or as a limit)
of filters.

 For the
sequential cascade $V$ we define the {\it contour} of $V$ (we write
$\int V$) inductively: if $r(V)=1$ then $\int V$ is a co-finite filter on $\max(V)$,
 if $W=V \looparrowleft V_w$ then $\int W = \int_{\int V} \int V_w$.
Similar
filters were considered in \cite{Kat1}, \cite{Kat2}, \cite{Dagu}.
Let $V$ be a monotone sequential cascade and let $u=\int V$.
Then
the {\it rank $r(\int V)$} of $\int V$
is, by definition, the rank of $V$.
It was shown in \cite{DolStaWat}, that for each countable ordinal $\alpha\geq 1$,
there is a monotone sequential contour of rank $\alpha$.
It
was shown in \cite{DolStaWat} that if $\int V= \int W$, then $r(V) =
r(W)$. The reader may find more
information about monotone sequential cascades and their contours in
\cite{Dol-multi}, \cite{DolMyn},
\cite{DolStaWat}, \cite{Star-ff}, \cite{Star-Ord-V-P},
\cite{Star-P-hier}.

We say that an ultrafilter $u$ belongs to a class $\cP_\alpha$ (we write $u\in \cP_\alpha$) if 

1) for each $\beta<\alpha$ there is a monotone sequential contour of rank $\beta$ contained in $u$

2) there is no monotone sequential contour of rank $\alpha$ contained in $u$.


Although this paper is self-contained we suggest to look at \cite{Star-Ord-V-P},
\cite{Star-P-hier} for more information concerning P-hierarchy.

IMPORTANT: In the remainder of this paper each filter is considered to be on
$\omega$, unless indicated otherwise, and for $f,g: \omega \rightarrow \omega$  we say that $f$ {\it dominates} $g$ if $f(n) > g(n)$ for all(!) $n<\omega$; this understanding does not change a domination number $\mathfrak{d}$.

%

\section{Results}

Let $V$ be a monotone sequential cascade of rank $\geq 2$. If from $V$ we remove all branches of height 1 obtaining a cascade $W$ then $\int V = \int W$. Thus we assume that each cascade of rank $\geq 2$ has no branches of height 1.

Let $V$ be a sequential cascade, we classically identify elements of the cascade with finite sequence of naturals by a function $f:V \rightarrow \omega^{<\omega}$ as follows:

$f(\emptyset _V) = \emptyset$; $f(w)=f(v)^\frown n$ if $w$ is the n-th element of $v^+$.  As a convention, we identify $v$ with $f(v)$, and see the cascade as a subset of $\omega^{<\omega}$.

A sequential cascade $V$ is {\it absorbing} if it fulfills the following condition:
if $(a_1,...,a_n)$ belongs to $V$ and $b_i\geq a_i$ for each $i\in \{1, ... , n\}$ then $(b_1,...,b_n) \in V$.
Note that each absorbing cascade is monotone. A contour of the absorbing cascade is called an absorbing contour.

\begin{rem}
For each monotone sequential cascade $V$ of rank less then, or equal to $\omega$
there is an absorbing cascade $W$ such that $\int V = \int W$. (It was proved in \cite{DolStaWat} that then $r(V)=r(W)$).
\end{rem}

\Proof: For $V$ of finite rank, it suffices to remove all branches of height less then $r(V)$ and re-enumerate branches.
Take a monotone sequential cascade $V$ of rank $\omega$, let $V = (n)  \looparrowleft V_n$ and in each $V_n$ remove all branches of height less then $r(V_n)$ to obtain the cascade we looking for (after re-enumerating of branches). $\Box$

We do not know weather we can extended Remark 2.1 to cascades of higher ranks $^1$ \footnotetext{1) We suppose that there is a counter-example for each $\alpha>\omega$}, but we have a little weaker Theorem 2.5 for them; first we need a lemma where by $-1+\gamma$ we denote $\gamma$ if $\gamma $ is infinite, $\gamma-1$ if $\gamma<\omega$.

\begin{lemm}
For each countable ordinal $\gamma$ there is a $(-1+\gamma +1)$-sequence $((a_{\alpha, \gamma}^n)_{n<\omega})_{1<\alpha\leq\gamma}$ of non decreasing $\omega$-sequences of ordinal numbers, such that 
$\lim_{n<\omega}(a_{\alpha, \gamma}^n+1)=\alpha$ and
$\alpha<\beta<\gamma$ implies
$a_{\alpha, \gamma}^n \leq a_{\beta, \gamma}^n$ for each natural number $n$.
\end{lemm}
\Proof: In the contrary. Let $\gamma$ be the first ordinal, such that the claimed sequence does not exist.

If $\gamma=\delta+1$ for some $\delta$, then for each $n$ it suffices to take :
$a_{\alpha, \gamma}^n=a_{\alpha, \delta}^n$ for each $\alpha \leq \delta$, and $a_{\gamma, \gamma}^n=\delta$.

If $\gamma$ is  a limit, take
an increasing sequence $(\gamma_n)$ of ordinals, such that $\gamma_1=1$ and  $\lim_{n<\omega}(\gamma_n+1)=\gamma$.
We define: 

$a_{\alpha, \gamma}^2=1$

for $\alpha$ such that $\gamma_n < \alpha \leq \gamma_{n+1}$, for $m\leq n$ let $a_{\alpha, \gamma}^m=\gamma_m$

for $\alpha$ such that $\gamma_n < \alpha \leq \gamma_{n+1}$, for $m>n$ let $a_{\alpha, \gamma}^m= \max\{\gamma_n, a^m_{\alpha, \gamma_{n+1}}\}$

Standard check shows that the defined sequence fulfills the claim. $\Box$

In forthcoming paper \cite{Star-Ord-V-P} we showed the following Remark 2.3 with standard proof by induction with respect to rank.

\begin{rem}\cite{Star-Ord-V-P}
For each monotone sequential cascade $V$, for each ordinal $\alpha < r(V)$
there exists a monotone sequential cascade $W$ of rank $\alpha$ such that $\int W \subset \int V$.
\end{rem}

\begin{lemm}[Folklore]
Let $(a_n)$ and $(b_n)$ be nondecreasing sequences of ordinals such that $a_1=b_1$ and $\lim_{n<\omega}a_n=\lim_{n<\omega}b_n$.
Then there is a non-decreasing, finite-to-one suriection $f:\omega\rightarrow\omega$ such that
$b_n\geq a_{f(n)}$.
\end{lemm}
\proof: Put $f(1)=1$,  now let $f(n)=f(n-1)+1$ if $a_{f(n-1)+1} \leq b_n$,
and $f(n)= f(n-1)$ otherwise. $\Box$

\begin{thm}
For each monotone sequential cascade $V$ there is a absorbing cascade $W$ such that $r(V)=r(W)$ and $\int W \subset \int V$.
\end{thm}

\Proof: Fix $\gamma$ and a $(-1+\gamma+1)$-sequence $((a_{\alpha, \gamma}^n)_{n<\omega})_{1<\alpha\leq\gamma}$ from the Lemma 2.2. We will show a little more, i.e. that for each monotone sequential cascade $V$ of rank $\gamma$ there is a absorbing cascade $W$ such that: $r(w^+_n)=a_{r(w),\gamma}^n$ for each $w \in W\setminus\max(W)$.

Since a set $\{v: v\in \emptyset_V^+, r(v)< a_{r(V),\gamma}^1\}$ is finite, so without loss of generality, we can assume that $r(V) \geq a_{r(V),\gamma}^1$ for all $v \in \emptyset_V^+$.

Put $b_n=r(\emptyset_V^{+,n})$ and $a_n = a_{\gamma,\gamma}^n$ and fix a function $f$ from the Lemma 2.4.

By Remark 2.3 for each $n<\omega$ there is a monotone sequential cascade $T_n$ such that $r(T_n)= a_{\gamma,\gamma}^n$
and $\int T_n \subset \int V_n$. 

Let $K_n$ be a cascade obtained from cascades $T_m$ for $m\in f^{-1}(n)$ by identifying all such $\emptyset_{T_m}$,
i.e. $\emptyset_{K_n}^+= \bigcup_{m\in f^{-1}(n)}\emptyset_{T_m}^+$, and if $k\in K$, $k\not= \emptyset_K$ then $k\in T_m$ for some $m\in f^{-1}(n)$ and $k^{+k}=k^{+T_m}$. Now, for each $K_n$ we use an inductive assumption, obtaining $W_n$, and the confluence $(k) \looparrowleft K_k$ is the cascade we are looking for. $\Box$

\begin{cor}
An ultrafilter $u$ contains a contour of a monotone sequential cascade of rank $\alpha$ if and only if $u$ contains a contour of an absorbing cascade of the rank $\alpha$.
\end{cor}

Let $V$ be a sequential cascade, and let $f,g: V\setminus \max(V) \rightarrow \omega$.
We say that a function $f$ $V$-{\it dominates} $g$
(in symbols $f \geq_V g$) if there is $U \in \int V$ such that
$f(v)\geq g(v)$ for each $v \in U^\downarrow$.
In this meaning we define the V-dominating family and analogically V-dominating number
${\mathfrak{d}}_V$. Define also the set $V(f)$ inductively:
$\emptyset_V\in V(f)$ and if $v \in V(f)$ then
$v^\frown k\in V(f)$ if $k\geq f(v)$. For $U \in \int V$ we define in the analogical way
$f_U: dom(f_U)\rightarrow \omega$: 
$\emptyset_V \in dom(f_U)$, if $v\in \dom(f_U)$ then
$f(v)=\min\{n<\omega: \forall m\geq n, \int (v^\frown n)^\uparrow \# U\}$;
 $v^\frown m \in \dom(f_U)$
for all $m\geq f_U(v)$.

\begin{rem}
For an absorbing cascade $V$ a family $\mathbb{F}$ of functions $V\setminus\max V \rightarrow \omega$ is V-dominating if and only if a family $\{ V(f): f \in \mathbb{F}\}$ is a base of $\int V$.
\end{rem}

For an absorbing cascade $V$ and for a function $f: V\setminus\max V \rightarrow \omega$, we define inductively a partial function $f\ss$ : $V \supset \dom(f\ss) \rightarrow V$ as follows:

$\emptyset_V \in \dom(f\ss)$ and $f\ss(\emptyset_V)=\emptyset_V$;

if $v \in \dom(f\ss)$ and $f\ss(v) \not\in \max V$ then 
$(v^\frown k)\in \dom(f\ss)$ for $k\geq f(v)$ and 
$f\ss(v^\frown k) =f\ss(v)^\frown (k-f(v)+1)$.

Note that $f\ss$ is a bijection and $\rng(f\ss)=V$.

If $g$ is also a function $V \setminus \max(V)\rightarrow \omega$, then $f\ss(g)(v):V\rightarrow \omega$ is defined by:
$f\ss(g)(v)=g(f\ss)^{-1}(v)$.

\begin{thm}
Let $V$ be an absorbing cascade, a family  $\mathbb{F} \subset \omega^{V \setminus \max(V)}$  is V-dominating if and only if $\mathbb{F}^*=\{f\ss(f):f\in F\}$ is a dominating family on $V \setminus \max(V)$.
\end{thm}

\Proof: Let $\mathbb{F} \subset ^{V \setminus \max(V)}\omega$ be a $V$-dominating family,
and let $f: V \setminus \max(V)\rightarrow \omega$. Take $g: V \setminus \max(V)\rightarrow \omega$ defined as follows:
$g(a_1, \ldots, a_n) = \max\{f((b_1, \ldots, b_n));$ $ b_i\leq a_i, (b_1, \ldots, b_n) \in V\}$. Since $\mathbb{F}$ is
$V$-dominating, there is $h\in \mathbb{F}$ that $V$-dominates $g$. Clearly $h\ss(h) \geq f$.

Now let $f$ be a witness that $\mathbb{F}^*$ is not dominating on $V$, define $g$ as above and observe that $g$ is not $V$-dominating by any $f \in \mathbb{F}$. $\Box$

By Corollary 2.6 Remark 2.7 and Theorem 2.8 we have: 
\begin{cor}
For each absorbing cascade V there is $\mathfrak{d}_V=\mathfrak{d}$.

\end{cor}

\begin{rem}[Folklore]
The minimum of cardinalities of families of non-dominating families such that the sum of all that families is a dominating family is $\mathfrak{b}$
\end{rem}
\Proof: For each $\alpha$, $\alpha<\lambda<\mathfrak{b}$, let $\mathbb{F}_\alpha$ be a non-dominating family of functions $\omega \rightarrow \omega$. Let $f_\alpha$ be a function non dominated by $\mathbb{F}_\alpha$, let $f$ be a function that dominates all $f_\alpha$ (there is some since $\lambda < \mathfrak{b}$). Clearly $f$ can not be dominated by any element of $\bigcup_{\alpha<\lambda}\mathbb{F}$.

Let $(f_\alpha)_{\alpha<\mathfrak{b}}$ be e non-limited sequence of functions. Define $\mathbb{F}_\alpha$ as a family
of all functions $f$ such that $f$ does not dominate $\{f_\beta: \beta \leq \alpha\}$. Clearly $\mathbb{F}_\alpha$
is non-dominating, and $\bigcup_{\alpha<\mathfrak{b}}\mathbb{F}_\alpha = \omega^\omega$ and so is dominating. $\Box$

\vspace{5mm}

Let $\mathbb{A} \subset ^X2$, a {\it supersets closure of }$\mathbb{A}$ is a family 
$\SC(\mathbb{A}) = \bigcup_{A\in \mathbb{A}}  \langle A \rangle$.
Let $u$ be a filter on $X$, we say that a family $\mathbb{P}\subset$ $2^X$ is a $\pi$-base of $u$
if $\mathbb{P}$ has finite intersection property, and if $u \subset \SC(\mathbb{P})$.

\begin{cor}
A sum of less then $\mathfrak{b}$ families $\mathbb{A}_\alpha$ that do not contain a $\pi$-base of absorbing contour of rank $\alpha>1$ does not contain any $\pi$-base of
monotone sequential contour of rank $\alpha$
\end{cor}
\Proof: Suppose on the contrary that there is a pair of witnesses -  a sequence $(\mathbb{A}_\alpha)_{\alpha<\lambda<\mathfrak{b}}$ 
and a monotone sequential cascade $W$, and fix a classical - defined by functions $V\setminus \max(V) \rightarrow \omega$ - base $\mathbb{B}$ of $V$.
By Theorem 2.5 or by Remark 2.1 (depending on the rank of the cascade) there is an absorbing cascade $V$ of the rank $r(V)=r(W)$ such that $\int V\subset \int W$, clearly each $\pi$-base of $\int W$ is also a $\pi$-base of $V$, so it suffices to prove this corollary for absorbing cascades. 
Since a family $\mathbb{A}_\alpha$ does not contain a $\pi$-base of $\int V$ and since $\mathbb{A}_\alpha$ is a $\pi$
-base of $\mathbb{B}_\alpha=\{B\in \mathbb{B}: B \supset A, A\in \mathbb{A}_\alpha\}$ , thus $\mathbb{B}_\alpha $ does not contain
$\int V$, and therefore by Remark 2.7 a family 
$\{f_B; B \in \mathbb{B}_\alpha\}$ is not V-dominating for each $\alpha<\lambda$. By Theorem 2.8
$\{f_B\ss(f_B); B \in  \mathbb{B}_\alpha \}$ is not dominating so by Remark 2.10 a family
$\bigcup_{\alpha<\lambda}\{f_B\ss(f_B); B \in  \mathbb{A}_\alpha \}$ is not dominating, and by Theorem 2.8 $\bigcup_{\alpha<\lambda}\{f_B; B \in  \mathbb{B}_\alpha \}$ is not V-dominating. Therefore by Remark 2.7 a family $\bigcup_{\alpha<\lambda} \mathbb{B}_\alpha $ does not contain a base of $\int V$, but since $\bigcup_{\alpha<\lambda}  \mathbb{B}_\alpha $ contains all supersets of elements of $\bigcup_{\alpha<\lambda} \mathbb{A}_\alpha$ which belong to  $\mathbb{B}$ thus 
$\bigcup_{\alpha<\lambda} \mathbb{B}_\alpha$ does not contain a $\pi$-base of $\int V$, and since $\bigcup_{\alpha<\lambda} \mathbb{A}_\alpha \subset \bigcup_{\alpha<\lambda} \mathbb{B}_\alpha$ thus
$\bigcup_{\alpha<\lambda} \mathbb{A}_\alpha$ does not contain a $\pi$-base of $\int V$.
$\Box$





\begin{cor}
An increasing ($\subset$) sequence of length less then $\mathfrak{b}$ of filters that do not contain a $\pi$-base of an absorbing sequential contour of rank $\alpha>1$ does not contain any 
monotone sequential contour of rank $\alpha$.
\end{cor}




\begin{rem} If $v$ is a filter and $T$ is a set such that $T\# v$ then if $v\mid_T$ does not contain any monotone sequential contour of rank $\alpha$ than $v$ does not contain any $\pi$-base of any monotone sequential contour of rank $\alpha$.
\end{rem}

\begin{rem}
Let $u$ be a filter which does not contain any $\pi$-base of  any monotone sequential cascade of rank $\alpha$, and let
 $f:\omega \rightarrow \omega$, then a filter
$\langle\{f^{-1}[U]: U\in u\}\rangle$
 does not contain any $\pi$-base of any monotone sequential cascade of rank $\alpha$.
\end{rem}

\Proof: For each $n$ such that $f^{-1}(n)$ is nonempty put $x_n = \min(f^{-1}(n))$ and let $X$ be a set of all such $x_n$'s. It suffices to consider $\langle\{f[U]: U\in u\}\rangle \mid_X$  which is a copy of $u$. $\Box$

\begin{thm}\cite[Theorem 2.5]{Star-P-hier} Let $(\alpha_n)_{n<\omega}$ be a non-decreasing
sequence of ordinals less than $\omega_1$, let $\alpha =
\lim_{n < \omega} (\alpha_n )$, let $1 < \beta <\omega_1$. If
$u_n \in \cP_{\alpha_n}$ is a discrete sequence of ultrafilters and
$u \in \cP_{\beta}$ then $\int_u u_n \in \cP_{\alpha+(-1+\beta)}$.
\end{thm}

\begin{rem}
Let $u$ be such a filter that there is a map $\omega \to \omega$ that $f(u)=\mathfrak{Fr}$. If $\langle u \cup \mathbb{A} \rangle$ is an ultrafilter for some family (of sets) $\mathbb{A}$ then $\card(\mathbb{A})\geq \mathfrak{u}$.
\end{rem}

\Proof: Let $f:\omega \to \omega$ be a function, such that $f(u)=\mathfrak{Fr}$, if $\langle u \cup \mathbb{A} \rangle$ is an ultrafilter thus $f(\langle u \cup \mathbb{A} \rangle)$ is a free ultrafilter and so $\{f[A]: A\in \mathbb{A}\}$ is a base of ultrafilter i.e. it has a cardinality of at least $\mathfrak{u}$. $\Box$.

\begin{thm} 
($\mathfrak{b}=\mathfrak{c}$) Let $1<\xi\leq\omega_1$ and let $p \in \cP_\xi$, then there is 
$\mathfrak{U}\subset \cP_\alpha$ of cardinality $\mathfrak{b}$ that 
$u >_{RK} p$  for each $u\in \mathfrak{U}$, and that elements of $\mathfrak{U}$ are Rudin-Kisler incomparable.
$^{1)}$ \footnotetext{1) We can obtain an easier version of the Theorem in a much shorter way:
(P-points exists) Let $\alpha$ be infinite countable ordinal, than for each $u\in\cP_\alpha$ there is $v\in \cP_\alpha$ such that $v>_{RK}u$.

\Proof: Let $\alpha\geq \omega$, and take any $u \in \cP_\alpha$.
 Consider a partition $(A_n)$ of $\omega$ into $\omega$ infinite sets.
 For each $n$ let $u_n$ be a P-point 
 such that $A_n \in u_n$. Put $v = \int_u u_n$, by Theorem 2.15 $v \in \cP_\alpha$,
 on the other hand for a function $f(m)=n$ for $m\in A_n$ we have $f(v)=u$ and there is no set
 $V\in v$ that $f\mid_V$ is one-to-one,  so it can not be Rudin-Kisler equivalent. 

For $\alpha = \omega_1$ the claim is obvious. $\Box$}

\end{thm}

\Proof: 

For $\xi=\omega_1$ the claim is obvious.
Fix $1<\xi< \omega_1$.
Let $f:\omega \rightarrow \omega$ be a finite-to-one function such that $\sup\{n:n \in U\}=\omega$ for each $U \in p$.

Let ${\cal{A}}_n$ be a family of such subsets of $\omega$ that there is $P\in p$ that $\card(f^{-1}(m))=\card(f^{-1}(m)\cap A)+n$ for each $m \in P$.

\begin{propo}
For each natural number $i$ a family $\mathbb{B}_i=\{f^{-1}[P] \cap A^c: P\in p, A \in {{\cal{A}}_i}\}$
does not contain a $\pi$-base of any monotone sequential contour of rank $\xi$. 
\end{propo}

On the contrary. Let $i=\min\{j<\omega:$ there is a $\pi$-base monotone sequential cascade of rank $\alpha$
contained in $\mathbb{B}_j\}$, let $\mathbb{P}$ be a $\pi$-base of absorbing sequential contour $\int V$ of rank $\xi$, such that $\mathbb{P}\subset \mathbb{B}_j$. For each $U \in \int V$, for each  $k\leq i$ let
$W_k(U)=\{n\in\omega: \card(f^{-1}(m))=\card(f^{-1}(m)\cap U)+k\}$, since $p$ is an ultrafilter, then for each $U\in \int V$ there is $k(U)\leq i$ such that $W_k(U)\in p$. Define also $A_k(U)=\{n<\omega: n\in f^{-1}(m)$ for such $m$ that $ \card(f^{-1}(m))=\card(f^{-1}(m)\cap U)+ k$ and $n \not\in U\}$. 
By minimality of $i$ there is $U_0 \in V$ such that $k(U_0)=i$, moreover if $U_1 \subset U_0$ then $k(U_1)=i$ and $A_i(U_1) \subset A_i(U_0)$.
Thus $\int V \subset \SC(\{f^{-1}[P] \cap A_i(U_0)^c: P\in p\})$ what is impossible since $\SC(\{f^{-1}[P] \cap A_i(U_0)^c: P\in p \}) = \SC(\{(f\mid_{A_i(U_0)^c})^{-1}[P]: P\in p\})$ so by Remark 2.14 it does not contain any monotone sequential cascade of rank $\xi$. $\Box$

\begin{propo}
$\langle f^{-1}[P] \cap A^c: P\in P, A\in \cal{A} \rangle$
 do not contain any monotone sequential contour of rank $\xi$.
\end{propo}

By Proposition 2.18 each $\mathbb{B}_i = \{f^{-1}[P] \cap A^c: P\in p, A \in {{\cal{A}}_i}\}$ do not contain a $\pi$-base of any monotone sequential contour of rank $\xi$,
thus by Corollary 2.11 $\bigcup_{n<\omega}\mathbb{B}_i$ does not contain a $\pi$-base of any monotone sequential contour of rank $\xi$,
but
$\bigcup_{n<\omega}\mathbb{B}_i = \langle f^{-1}[P] \cap A^c: P\in P, A\in \cal{A} \rangle$,
i.e., is a filter and so since it does not contain a $\pi$-base of any monotone sequential contour of rank $\xi$, it also does not contain any monotone sequential contour of rank $\xi$. $\Box$

Clearly there is no set $\bar{A}$ that $\langle f^{-1}[P] \cap A^c \cap \bar{A}: P\in P, A\in \cal{A} \rangle$ is an ultrafilter thus there is a sequence of pairwise disjoint sets $(C_n)_{n<\omega}$ such that $C_n\# 
\langle f^{-1}[P] \cap A^c: P\in P, A\in {\cal{A}} \rangle$ and $\bigcup_{n<\omega}C_n  = \omega$. Let ${\cal T} = \langle \{f^{-1}[P] \cap A^c: P\in P, A\in {\cal{A}}\}\cup \{\bigcup_{n>m} C_n: m<\omega\}\rangle$.

\begin{rem}
A filter $\langle {\cal T} \cup \bar{A}\rangle$ 
dose not contain any monotone sequential contour of rank $\xi$ for any $\bar{A}\#{\cal T}$, and there is a function $f:\omega \to \omega$ that $f({\cal T})=\mathfrak{Fr}$.
\end{rem}
First part follows from Proposition 2.19, and a number of generators added to the family $\langle f^{-1}[P] \cap A^c \cap \bar{A}: P\in P, A\in \cal{A} \rangle$ in the virtue of Corollary 2.11. To see the second part of Remark 2.20 it suffices to consider a function $f:\omega \to \omega$ such that $f(n)=m$ if $n\in C_m$. $\Box$

We enlist all absorbing cascades of rank $\xi$ in a sequence $(V_\alpha)_{\alpha<\mathfrak{b}}$ and all functions $\omega \to \omega$ in a sequence $(f_\beta)_{\beta<\mathfrak{b}}$.
We will build a family $\{({\cal{F}}_\alpha)^\beta\}_{\beta<\mathfrak{b}}$
of increasing $\mathfrak{b}$-sequences $({\cal{F}}_\alpha)_{\alpha<\mathfrak{b}}$ of filters such that:

1)Each ${\cal{F}}_\alpha^\beta$ is generated by $\cal T$ together with some family of cardinality $< \mathfrak{b}$ of sets;

2) ${\cal{F}}^\beta_0 = {\cal T}$ for each $\beta<\mathfrak{b}$;

3) For each $\alpha, \beta < \mathfrak{b}$, there is $F\in {\cal{F}}_{\alpha+1}^\beta$ such that $F^c\in \int V_\alpha$;

4) For limit $\alpha$ for each $\beta$, ${\cal{F}}_\alpha^\beta=\bigcup_{\gamma<\alpha}{\cal{F}}_\gamma^\beta$;

5) For each $\alpha$, for each $\gamma<\alpha$, for each $\beta_1$, $\beta_2 <\alpha$ 
there is a set $F\in {\cal{F}}_{\alpha+1}^{\beta_1}$,  such that $(f_\gamma[F])^c\in {\cal{F}}_{\alpha+1}^{\beta_2}$.

Existence of such families is a standard work by induction with respect to $\alpha$ with sub-induction with respect to $\gamma<\alpha$, with sub-sub induction with respect to $\beta_1<\gamma$ and with sub-sub-sub-induction with respect to $\beta_2<\beta_1$, using Remark 2.16, Remark 2.20 and Remark 2.10.

Now it suffice for each $\beta < \mathfrak{b}$ take any ultrafilter extending $\bigcup_{\beta<\mathfrak{b}}{\cal{F}}_\alpha^\beta$. $\Box$

 \begin{thm}
($\mathfrak{b}=\mathfrak{c}$) Let $1<\xi\leq \omega_1$, and let $(p_n)_{n<\omega}$ be a RK-increasing sequence of elements of $\cP_\xi$, then there exists $u\in \cP_\xi$ such that $u >_{RK} p_n$ for each 
 $n<\omega$.
\end{thm}

\Proof: Let $f_n$ be a function $\omega \rightarrow \omega$ - witness that $p_{n+1} >_{RK} p_n$.
For each natural number $m$ consider on $\omega \times \omega$ a family of sets ${\mathbb{B}}_m$ 
such that 
${\mathbb{B}}_m\mid (\omega \times \{n\})$
$= \langle \{f_{n-1}^{-1}\circ f_{n-1}^{-1} \circ \ldots \circ f_m^{-1}(P): P\in p_n\}\rangle$ for $n\geq m$.
Let $\mathbb{B}=\bigcup_{m<\omega} \mathbb{B}_m$. Clearly $\mathbb{B}$ is a filter, and each ultrafilter which extends 
$\mathbb{B}$ is RK greater then each $p_n$. Also, by Remark 2.14, each $\mathbb{B}_n$ does not contain any monotone sequential contour of rank $\xi$ so by Corollary 2.12, $\mathbb{B}$ does not contain any monotone sequential contour of rank $\xi$.

We enlist all absorbing cascades of rank $\xi$ in a sequence $(V_\alpha)_{\alpha<\mathfrak{b}}$ and we will build an increasing $\mathfrak{b}$-sequence of filters ${\cal{F}}_\alpha$ such that:

1) ${\cal{F}}_0= \mathbb{B}$.

2) For each $\alpha$, there is such $F\in {\cal{F}}_{\alpha+1}$ that $F^c\in \int V_\alpha$;

3) For a limit $\alpha$, ${\cal{F}}_\alpha=\bigcup_{\beta<\alpha}{\cal{F}}_\beta$.

The rest of the proof is an easier version of the final part of the proof of Theorem 2.17. $\Box$

\begin{cor}
($\mathfrak{b}=\mathfrak{c}$) Let $1<\xi\leq \omega_1$, and let $(p_n)_{n<\omega}$ be a RK-increasing sequence of elements of $\cP_\xi$, then there exists a family 
$\mathfrak{U}\subset \cP_\xi$ of cardinality $\mathfrak{b}$
such that $u >_{RK} p_n$ for each 
 $u\in \mathfrak{U}$, and that elements of $\mathfrak{U}$ are Rudin-Kisler incomparable.
\end{cor}
\Proof: Just combine Theorem 2.21 with Theorem 2.17. $\Box$

\begin{thm}\cite[Theorem 2.8]{Star-P-hier}
The following statements are equivalent:

1) P-points exist;

2) Classes $\cP_\alpha$ are nonempty for each countable successor
$\alpha$;

3) There exists  a countable  successor $\alpha>1$ such that the
class $\cP_\alpha$ is nonempty.

\end{thm}

\begin{thm}\cite{Ket}
$\mathfrak{d}=\mathfrak{c}$ if, and only if, every filter generated by less then $\mathfrak{c}$ elements can be extended to a
P-point.
\end{thm}

\begin{thm}
($\mathfrak{b}=\mathfrak{c}$) Each class of $\cP$-hierarchy is non-empty.
\end{thm}
\Proof:
For classes of index $\xi \in \{1,\omega_1\}$ we are done (in ZFC) by \cite[Corollary 6.4]{Star-Ord-V-P}

For successor $1<\xi<\omega_1$, since $\mathfrak{b}\leq \mathfrak{d}$, it suffice to combine Theorem 2.23 and Theorem 2.24.

For limit $\xi<\omega_1$ a proof is essentially the same as final part of the proof of Theorem 2.17:

We enlist all absorbing cascades in a sequence $(V_\alpha)_{\alpha<\mathfrak{b}}$.
Let $(\int V_n)$ be an increasing ("$\subset$") sequence of monotone sequential contours such that
$\lim_{n <\omega}(r(V_n)+1)=\xi$, such sequences exist in ZFC - and were constructed in the proof of \cite[Theorem 4.6]{DolStaWat}. Thus by Corollary 2.11 $\bigcup_{n<\omega}\int V_n$ does not contain a $\pi$-base of monotone sequential cascade of rank $\xi$.
 
We enlist all absorbing cascades of rank $\xi$ in a sequence $(V_\alpha)_{\alpha<\mathfrak{b}}$ and we will build an increasing $\mathfrak{b}$-sequence of filters ${\cal{F}}_\alpha$ such that:

1) ${\cal{F}}_0= \bigcup_{n<\omega}\int V_n$.

2) For each $\alpha$, there is $F\in {\cal{F}}_{\alpha+1}$ such that $F^c\in \int V_\alpha$;

3) For a limit $\alpha$, ${\cal{F}}_\alpha=\bigcup_{\beta<\alpha}{\cal{F}}_\beta$.

The rest of the proof is an easier version of the final part of the proof of Theorem 2.17. $\Box$

\begin{cor} 
($\mathfrak{b}=\mathfrak{c}$)
Each class of P-hierarchy of index $>1$ contains a family of cardinality $\mathfrak{b}$ of pairwise Rudin-Kisler incomparable ultrafilters. 
\end{cor}

\Proof: Just combine Theorem 2.25 with Theorem 2.17. $\Box$

 \hspace{-2mm}
\gd

\medskip

{\small\sc \noindent Andrzej Starosolski, {Wydzia\l}  {Matematyki Stosowanej},
Politechnika \'{S}l\c{a}ska, Gliwice, Poland

E-mail:  andrzej.starosolski@polsl.pl}

\end{document}